\documentclass[11pt]{amsart}
\usepackage{amsxtra}
\usepackage{amssymb}
\theoremstyle{plain}

\newcommand{\cleqn}{\setcounter{equation}{0}}
\newcommand{\clth}{\setcounter{theorem}{0}}
\newcommand {\sectionnew}[1]{\section{#1}\cleqn\clth}

\newtheorem{theorem}{Theorem}[section]
\newtheorem{lemma}[theorem]{Lemma}
\newtheorem{definition-theorem}[theorem]{Definition-Theorem}
\newtheorem{proposition}[theorem]{Proposition}
\newtheorem{corollary}[theorem]{Corollary}
\newtheorem{definition}[theorem]{Definition}
\newtheorem{example}[theorem]{Example}
\newtheorem{remark}[theorem]{Remark}
\newtheorem{conjecture}[theorem]{Conjecture}
\newtheorem{notation}[theorem]{Notation}
\newcommand \bth[1] { \begin{theorem}\label{t#1} }
\newcommand \ble[1] { \begin{lemma}\label{l#1} }

\newcommand \bpr[1] { \begin{proposition}\label{p#1} }
\newcommand \bco[1] { \begin{corollary}\label{c#1} }
\newcommand \bde[1] { \begin{definition}\label{d#1}\rm }
\newcommand \bex[1] { \begin{example}\label{e#1}\rm }
\newcommand \bre[1] { \begin{remark}\label{r#1}\rm }
\newcommand \bcj[1] { \begin{conjecture}\label{j#1}\rm }

\newcommand \bnota[1] { \begin{notation}\label{n#1}\rm }
\renewcommand {\eth} { \end{theorem} }
\newcommand {\ele} { \end{lemma} }

\newcommand {\epr} { \end{proposition} }
\newcommand {\eco} { \end{corollary} }
\newcommand {\ede} { \end{definition} }
\newcommand {\eex} { \end{example} }
\newcommand {\ere} { \end{remark} }
\newcommand {\ecj} { \end{conjecture} }

\newcommand {\enota} { \end{notation} }
\newcommand \thref[1]{Theorem \ref{t#1}}
\newcommand \leref[1]{Lemma \ref{l#1}}

\newcommand \lb[1]{\label{#1}}

\def \Cset {{\mathbb C}}
\def \KK {{\mathbb K}}
\def \Zset {{\mathbb Z}}

\def \Qset {{\mathbb Q}}

\def \AA  {{\mathcal{A}}}           
\def \BB  {{\mathcal{B}}}
\def \EE  {{\mathcal{E}}}

\def \OO {{\mathcal{O}}}

\def \UU {{\mathcal{U}}}
\def \RR {{\mathcal{R}}}

\def \De {\Delta}   
\def \de {\delta}
\def \al {\alpha}

\def \la {\lambda}

\def \om {\omega}
\def \ga {\gamma}
\def \de {\delta}
\def \Ga {\Gamma}

\def \lha {\leftharpoonup}
\def \rha {\rightharpoonup}

\def \ci  {\circ}

\def \rcor {\rangle}
\def \lcor {\langle}

\def \ol {\overline}

\def \id { {\mathrm{id}} }

\def \Lie { {\mathrm{Lie \,}} }

\def \g  {\mathfrak{g}}   

\def \n  {\mathfrak{n}}

\DeclareMathOperator \Span { {\mathrm{Span}} }

\DeclareMathOperator \res { {\mathrm{res}} }

\newcommand \Spec { {\mathrm{Spec}} }
\begin{document}
\title[Strata of prime ideals]
{Strata of prime ideals of De Concini--Kac--Procesi algebras and 
Poisson geometry} 
\author[Milen Yakimov]{Milen Yakimov}
\dedicatory{To Ken Goodearl on his 65th birthday}
\address{
Department of Mathematics \\
Louisiana State Univerity \\
Baton Rouge, LA 70803 and
Department of Mathematics \\
University of California \\
Santa Barbara, CA 93106 \\
U.S.A.
}
\email{yakimov@math.lsu.edu}
\date{}
\begin{abstract}
To each simple Lie algebra $\g$ and an element $w$ of the corresponding
Weyl group De Concini, Kac and Procesi associated a subalgebra $\UU^w_-$
of the quantized universal enveloping algebra $\UU_q(\g)$, which is 
a deformation of the universal enveloping algebra $\UU(\n_- \cap w(\n_+))$
and a quantization of the coordinate ring of the Schubert cell corresponding
to $w$. The torus invariant prime ideals of these algebras were classified 
by M\'eriaux and Cauchon \cite{MC}, and the author \cite{Y}. 
These ideals were also explicitly described 
in \cite{Y}. They index the the Goodearl--Letzter strata of the stratification 
of the spectra of $\UU^w_-$ into tori. In this paper we derive a formula 
for the dimensions of these strata and the transcendence 
degree of the field of rational Casimirs on any open Richardson
variety with respect to the standard Poisson structure \cite{GY}.
\end{abstract}
\maketitle
\sectionnew{Introduction}
\lb{intro}
Assume that $A$ is a Noetherian $\KK$-algebra equipped with a rational action of 
a $\KK$-torus $T$ by algebra automorphisms. Under very general 
assumptions Goodearl and Letzter \cite{GL} constructed 
a stratification of $\Spec A$ into tori indexed by the $T$-primes of $A$. 
Previously Joseph \cite{J} and Hodges--Levasseur--Toro \cite{HLT} 
obtained such stratifications of the spectra of 
quantized coordinate rings of simple groups. The Goodearl--Letzter 
results showed that such stratifications of spectra of rings 
exist in much greater generality, in particular whenever 
$A$ is an iterated skew polynomial extension under some natural 
assumptions relating the structure of $A$ and the action of $T$.
This generated a lot of research in ring theory targeted at 
the explicit description of the above stratification of $\Spec R$
for concrete rings $R$. The Cauchon approach of deleted derivations 
\cite{C} provides an iterative procedure to classify the $T$-primes 
of an iterated skew polynomial extension. The explicit realization 
of this procedure often leads to difficult combinatorial problems. 
After many specific rings were investigated, most notably
the algebras of quantum matrices (see \cite{C,GLe,La}), 
it was observed that almost all of 
them fit into the general class of quantized universal enveloping 
algebras of nilpotent (or slightly more generally solvable) Lie algebras. 

The most general quantization of a class of nilpotent Lie algebras 
up to date was constructed by De Concini, Kac and Procesi \cite{DKP} using the 
Lusztig root vectors of a universal enveloping algebra $\UU_q(\g)$ 
of a simple Lie algebra $\g$. The algebras are parametrized 
by elements $w$ of the Weyl group $W_\g$ of $\g$. The corresponding 
algebra $\UU^w_-$ is a deformation 
of the universal enveloping algebra of $\n_- \cap w(\n_+)$ 
where $\n_\pm$ is a pair of opposite nilpotent subalgebras of 
$\g$. It can be viewed \cite{DP,Y} as a quantization 
of the coordinate ring of the Schubert cell 
$B_+ w \cdot B_+$ in the full flag variety $G/B_+$ 
with respect to the standard Poisson structure \cite{GY}.
Within the framework of quantum flag varieties,
it was proved that the quantum Schubert cells 
of all quantum partial flag varieties 
are isomorphic to algebras of the form $\UU^w_-$, 
see \cite[Theorem 3.6]{Y2}, and this was used to 
classify the torus invariant prime ideals of all quantum 
flag varieties \cite{Y2}. 

Denote by $T$ the maximal torus of the 
connected simply connected algebraic group corresponding to $\g$.
In \cite{MC} M\'eriaux and Cauchon classified the $T$-primes of 
$\UU^w_-$ for $q \in \KK^*$ not a root of unity and arbitrary base field
$\KK$. There is no explicit description of the $T$-primes 
in this approach, in particular the inclusions between 
the $T$-primes in the M\'eriaux--Cauchon picture are currently 
unknown.

In \cite{Y} the author constructed a second realization of 
the algebras $\UU^w_-$ and derived an explicit description of 
each $T$-prime of $\UU^w_-$, based on works of Gorelik and Joseph
\cite{G,J0,J}. This resulted in a new parametrization 
of the set of $T$-primes of $\UU^w_-$ which also explicitly identified
the poset structure of the $T$-spectrum with the inclusion relation 
with a particular Bruhat interval in the related Weyl group. The statement
of these results are summarized in \cite[Theorem 1.1]{Y}. Although
these results were stated in \cite{Y} for fields $\KK$ of 
characteristic $0$ and deformation parameters $q$ which are 
transcendetal over $\Qset$, the proofs work in the general 
case without restrictions on $\KK$ assuming only that $q$ 
is not a root of unity. This will be addressed in a 
forthcoming preprint. The $T$-primes of $\UU^w_-$ are indexed 
by the elements $y \in W_\g$ such that $y \leq w$ under the
Bruhat order in $W_\g$. We denote by $I_w(y)$ the $T$-prime
ideal of $\UU^w_-$ corresponding to $y$ as in 
\cite[Theorem 1.1]{Y}, and refer the reader to 
\S \ref{qalg} for the explicit description of $I_w(y)$.

In this paper we derive a formula for the dimension
of the Goodearl--Letzter stratum of $\UU^w_-$ corresponding 
to each $T$-prime $I_w(y)$. This is done under the assumption 
that the base field $\KK$ has characteristic $0$ and the 
deformation parameter $q$ is transcendental over $\Qset$, 
see \thref{main}. The proof uses Poisson geometry.

Bell, Casteels, and Launois \cite{BCL} obtained 
simultaneously and independently a formula 
for the dimensions of the Goodearl--Letzter stratum of $\UU^w_-$
in the M\'eriaux--Cauchon parametrization \cite{MC}
for an arbitrary base field $\KK$, and  
$q \in \KK^*$ is not a root of unity. They gave a second proof 
for the case of quantum matrices in \cite{BCL2}. 
Previously the dimension of the stratum for the 
special case of the $0$ ideal of $\UU^w_-$ 
was obtained by Bell and Launois in \cite{BL}. 
The two approaches of M\'eriaux and Cauchon \cite{Y}
and the author \cite{Y} to the $T$-spectrum of $\UU^w_-$ 
are of very different nature and are not connected yet.
In particular, one cannot transfer the explicit 
description of $T$-primes and their inclusions from 
the picture in \cite{Y} to the picture in \cite{MC}.
In a forthcoming preprint 
we will describe a direct ring theoretic derivation of the 
dimension formulas for the Goodearl--Letzter strata 
of $\UU^w_-$ in the picture in \cite{Y} which works 
in the more general case of an arbitrary base field 
$\KK$ and $q$ which is not a root of unity. This can 
be also obtained by combining the results of this paper 
and \cite{BCL}, though leading to an indirect proof.

To describe in more concrete terms the results and methods of the paper
we return to the general setting of a $\KK$-algebra $A$ with an action of a 
torus $T$ by algebra automorphisms. Denote by $T-\Spec A$ the set 
of $T$-primes of $A$. For a $T$-prime $I$ denote 
the Goodearl--Letzter stratum \cite{GL} of $\Spec A$: 
\[
\Spec_I A = \{ J \in \Spec A \mid \cap_{t \in T} t \cdot J = I \}. 
\] 
Following \cite{GL}, consider the localization 
$(A/I)[\EE_I^{-1}]$ by the set of all homogeneous nonzero elements
$\EE_I \subset A/I$. Goodearl and Letzter proved (under certain mild assumptions) 
that 
\[
\Spec A = \bigsqcup_{I \in H-\Spec} \Spec_I A,
\]
$Z((A/I)[\EE_I^{-1}])$ is isomorphic to a Laurent polynomial ring, 
and that $\Spec_I A$ is homeomorphic to the torus 
$\Spec Z((A/I)[\EE_I^{-1}])$, see \cite[Theorem 6.6]{G0}.
This procedure was further developed 
in the book of Brown and Goodearl \cite[Part II]{BG} and by Goodearl 
\cite{G0} who also showed that the same holds true for smaller
sets $\EE_I$ which match the original 
Joseph \cite{J0} and Hodges--Levasseur--Toro \cite{HLT} methods,
see \cite[Theorem 5.3]{G0}. This partition has all topological 
properties required for a stratification, see \cite[Lemma 3.4]{G0}.
Brown and Goodearl also proved \cite[Theorem II.6.4]{BG} that, if 
$A$ is an iterated skew polynomial algebra satisfying certain 
general conditions which hold for the algebras $\UU^w_-$, then 
the base field for all Laurent polynomial rings 
$Z((A/I)[\EE_I^{-1}])$ is $\KK$. Since all algebras 
$\UU^w_\pm$ are iterated skew polynomial algebras, 
this fact applies to all of them.

In \cite{Y} we realized the algebras $\UU^w_-$ as explicit
quotients of Joseph's quantum translated Bruhat cell algebras, see
\cite[\S 10.4.8]{J} and \cite{G}. The latter 
are defined in terms of the quantized coordinate ring
$\OO_q(G)$ and satisfy similar commutation relations derived from 
the quantum $R$-matrix for $\g$. In this paper,
starting from this realization 
we construct enough elements in the center 
$Z((\UU^w_-/I_w(y))[\EE_{I_w(y)}^{-1}])$ and prove that this center 
is a Laurent polynomial ring in at least $\dim E_{-1} (w^{-1} y)$
variables. Here and below for a linear operator $L$ acting on a 
vector space $V$ and $c \in \Cset$ we denote by 
$E_c(L)$ the $c$-eigenspace of $L$. Next we pass to 
an integral form of $\UU^w_-/I_w(y)$ and show that its 
specialization at $q=1$ is isomorphic to the coordinate ring 
of the open Richardson variety 
\[
R_{y,w} = B^- y \cdot B^+ \cap B^+ w \cdot B^+ \subset G/B^+
\]
with Poisson algebra structure coming from the standard 
Poisson structure on $G/B^+$, see \cite{GY} and \S \ref{Poisson}.
If the dimension of the stratum $\Spec_{I_w(y)} \UU^w_-$
was strictly greater than $\dim E_{-1} (w^{-1} y)$, then 
the transcendence degree of the center 
of the Poisson field of rational functions on $R_{y,w}$ 
would be strictly greater than $\dim E_{-1} (w^{-1} y)$ which
is shown to contradict with the dimension formulas for
the symplectic leaves in $R_{y,w}$. This implies that both
the dimension of the Goodearl--Letzter stratum 
$\Spec_{I_w(y)} \UU^w_-$ and the transcendence degree of 
the center of the Poisson field of rational functions on 
the open Richardson variety $R_{y,w}$ are equal to  
$\dim E_{-1} (w^{-1} y)$. 
 
In \cite{Y3} we prove some further properties of the 
algebras $\UU^w_-$. Firstly, for each torus invariant prime ideal 
$I_w(y)$ we construct efficient polynormal generating sets.
In the special case of the algebras of quantum matrices
this leads to an explicit proof of the Goodearl--Lenagan 
conjecture \cite{GLe} that all torus invariant prime 
ideals of the algebras of quantum matrices have polynormal 
generating sets consisting of quantum minors. Furthermore,
we prove that all $\Spec \UU^w_-$ are normally separated,
and that each algebra $\UU^w_-$ is catenary.
 
The paper is organized as follows. Sect. \ref{qalg} 
contains background for quantized universal enveloping 
algebras, the algebras $\UU^w_-$ and their spectra.
Sect. \ref{Poisson} deals with the related Poisson 
structures on flag varieties. Sect. \ref{dimensions} 
carries out the connection between the two and contains  
the proofs of the main results.
\\ \hfill \\
{\bf Acknowledgements.} I am grateful to Ken Goodearl
who greatly inspired my interest in ring theory and 
taught me ring theory. I would also like to thank the 
referee for his/her thorough reading of the paper and his
questions which helped me to improve the exposition.

The research of the author was 
supported by NSF grants DMS-0701107 and DMS-1001632.
\sectionnew{Quantized nilpotent algebras}
\lb{qalg}
Fix a base field $\KK$ of characteristic $0$ 
and $q \in \KK$ which is transcendetal over $\Qset$. 
Let $\g$ be simple Lie algebra 
of rank $r$ and Cartan matrix $c_{ij}$. Denote by $\UU_q(\g)$ 
the quantized universal enveloping algebra of $\g$ over 
$\KK$ with deformation parameter $q$. 
It is a Hopf algebra over $\KK$ with generators
\[
X^\pm_i, K_i^{\pm 1}, \; i=1, \ldots, r,
\]
and relations
\begin{align*}
&K_i^{-1} K_i = K_i K^{-1}_i = 1, \, K_i K_j = K_j K_i, \,
K_i X_j^\pm K^{-1}_i = q_i^{\pm c_{ij}} X_j^\pm,
\\
&X^+_i X^-_j - X^-_j X^+_i = \de_{i,j} \frac{K_i - K^{-1}_i}
{q_i - q_i^{-1}},
\\
&\sum_{k=0}^{1-c_{ij}}
\begin{bmatrix}
1-c_{ij} \\ k
\end{bmatrix}_{q_i}
      (X^\pm_i)^k X^\pm_j (X^\pm_i)^{1-c_{ij}-k} = 0, \, i \neq j.
\end{align*}
where $d_i$ are positive relative prime integers such that 
$(d_i c_{ij})$ is a symmetric matrix and $q_i = q^{d_i}$. 
Denote by $\UU_\pm$ the subalgebras of $\UU_q(\g)$
generated by $\{X^\pm_i\}_{i=1}^r$. Let $H$ be the group 
generated by $\{K_i^{\pm 1}\}_{i=1}^r$ (i.e. the group 
of group like elements of $\UU_q(\g)$).

Let $P$ and $P_+$ be the sets of integral and dominant integral 
weights of $\g$. The sets of simple roots, simple coroots, and 
fundamental weights of $\g$ will be denoted by $\{\al_i\}_{i=1}^r,$ 
$\{\al_i\spcheck\}_{i=1}^r,$ and $\{\om_i\}_{i=1}^r,$ respectively.
Let $\lcor.,. \rcor$ be the nondegenerate invariant bilinear form
on $\g$ such that the square length of a short root is equal to 2.

The $q$-weight spaces of an $H$-module $V$ are defined by
\[
V_\la = \{ v \in V \mid K_i v = q^{ \lcor \la, \al_i \rcor} v, \; 
\forall i = 1, \ldots, r \}, \; \la \in P.
\]
A $\UU_q(\g)$-module is called a type 1 module if it is the 
sum of its $q$-weight spaces. Recall that all finite dimensional 
type 1 $\UU_q(\g)$-modules are completely reducible and are 
parametrized by $P_+$, see \cite[\S 10.1]{CP}. 
Let $V(\la)$ denote the irreducible 
weight $\UU_q(\g)$ module of highest weight $\la \in P_+$. 
For each $\la \in P_+$ fix a highest weight vector $v_\la$ 
of $V(\la)$ such that $\forall \la, \mu \in P_+$
$v_{\la+\mu} = v_\la \otimes v_\mu$ when $V(\la+\mu)$ is 
realized as a submodule of $V(\la) \otimes V(\mu)$.
 
All duals of finite dimensional $\UU_q(\g)$ modules
will be thought of as left modules using the antipode of $\UU_q(\g)$.

Denote the Weyl and braid groups of $\g$ by $W_\g$ and $\BB_\g$, 
respectively. Let $s_1, \ldots, s_r$ be the simple reflections 
of $W_\g$ corresponding to the roots $\al_1, \ldots, \al_r$, and
$T_1, \ldots, T_r$ be the standard generators of $\BB_\g$. 
Recall that one has a section $W_\g \to \BB_\g$, $w \mapsto T_w$,
of the canonical projection $\BB_\g \to W_\g$. Given a 
a reduced expression $w = s_{i_1} \ldots s_{i_k}$ one sets
$T_w = T_{i_1} \ldots T_{i_k}$. The latter does not depend 
of the choice of a reduced expression of $w$.

There are natural actions of $\BB_\g$ on $\UU_q(\g)$ 
and the modules $V(\la)$, see \cite[\S 5.2 and \S 37.1]{L} 
for details. They have the properties that 
$T_w ( x . v ) = (T_w x) . (T_w v)$, 
$T_w(V(\la)_\mu) = V(\la)_{w \mu}$
for all $w \in W_\g$, $x \in \UU_q(\g)$, 
$\la \in P_+$, $v \in V(\la)$, $\mu \in P$.

For a reduced decomposition
\begin{equation}
\label{wdecomp}
w = s_{i_1} \ldots s_{i_k}
\end{equation}
of an element $w \in W_\g$, define the roots
\begin{equation}
\label{beta}
\beta_1 = \al_{i_1}, \beta_2 = s_{i_1} \al_{i_2}, 
\ldots, \beta_k = s_{i_1} \ldots s_{i_{k-1}} \al_{i_k}
\end{equation}
and Lusztig's root vectors
\begin{equation}
X^{\pm}_{\beta_1} = X^{\pm}_{i_1}, 
X^{\pm}_{\beta_2} = T_{s_{i_1}} X^\pm_{i_2}, 
\ldots, X^\pm_{\beta_k} = T_{s_{i_1} \ldots s_{i_{k-1}}} X^{\pm}_{i_k},
\label{rootv}
\end{equation}
see \cite[\S 39.3]{L}. De Concini, Kac and Procesi defined 
\cite{DKP} the subalgebras $\UU_\pm^w$ of $\UU_\pm$ generated by 
$X^{\pm}_{\beta_j}$, $j=1, \ldots, k$ and proved:

\bth{DKP} (De Concini, Kac, Procesi) \cite[Proposition 2.2]{DKP}
The algebras $\UU_\pm^w$ do not depend on the choice of a 
reduced decomposition of $w$ and have the PBW basis
\begin{equation}
\label{vect}
(X^\pm_{\beta_k})^{n_k} \ldots (X^\pm_{\beta_1})^{n_1}, \; \; 
n_1, \ldots, n_k \in \Zset_{\geq 0}.
\end{equation}
\eth

The fact that the space spanned by the monomials 
\eqref{vect} does not depend on the choice of a
reduced decomposition of $w$ was independently obtained by Lusztig 
\cite[Proposition 40.2.1]{L}.

The quantized coordinate ring $R_q[G]$ of the split, connected, simply 
connected algebraic group $G$ corresponding to $\g$ is the Hopf  
subalgebra of the restricted dual of $\UU_q(\g)$ spanned by 
all matrix entries $c_{\xi, v}^\la,$ 
$\la \in P_+$, $v \in V(\la), \xi \in V(\la)^*$:
$c_{\xi,v}^\la(x) = \lcor \xi, x v \rcor$ for $x \in \UU_q(\g)$.
Denote by $R^+$ the subalgebra of $R_q[G]$ spanned by 
all matrix entries $c_{\xi , v_\la}^\la$ where $\la \in P_+$, 
$\xi \in V(\la)^*$ and $v_\la$ is the fixed highest weight vector 
of $V(\la)$. The group $H$ acts on $R_q[G]$ on the left 
and right by
\begin{equation}
\label{Hact}
x \rha c = \sum c_{(2)}(x)c_{(1)}, \; 
c \lha x = \sum c_{(1)}(x)c_{(2)}, \; 
x \in H, c \in R_q[G]
\end{equation} 
in terms of the standard notation for the comultiplication 
$\Delta(c) = \sum c_{(1)} \otimes c_{(2)}.$

For all $\la \in P_+$ and $w \in W_\g$ define 
$\xi_{w, \la} \in (V(\la)^*)_{- w\la}$ such that
$\lcor \xi_{w, \la}, T_w v_\la \rcor =1$. Let
\begin{equation}
\label{cwla}
c^\la_w = c^\la_{\xi_{w,\la}, v_\la}.
\end{equation}
Then $c^\la_w c^\mu_w = c^{\la+ \mu}_w= c^\mu_w c^\la_w$,
$\forall \la, \mu \in P_+$. Set
\[
c_w = \{ c_w^\la \mid \la \in P_+ \}. 
\]
It is \cite[Lemma 9.1.10]{J} an Ore subset of $R^+$. 

Denote the localization
\[
R^w = R^+[c_w^{-1}]
\]
(see \cite{J, G}) and observe that \eqref{Hact} induces
$H$-actions on $R^w$. The invariant subalgebra with respect 
to the left action \eqref{Hact} of $H$ is denoted by $R_0^w$ and is called 
the quantum translated Bruhat cell \cite[\S 10.4.8]{J}. It was studied in great
detail by Gorelik in \cite{G}.
For $\la \in P_+$ set $c_w^{-\la} = (c_w^\la)^{-1} \in R^w$. 
We have that
\begin{equation}
\label{Rw0}
R^w_0 = \{ c_w^{-\la} c^\la_{\xi, v_\la} \mid 
\la \in P_+, \xi \in V(\la)^* \}
\end{equation}
since 
\begin{equation}
\label{span}
\forall \la, \mu \in P_+, \; 
\xi \in V(\la)^*, \quad \exists \xi' \in V(\la+\mu)^*
\; \; 
\mbox{such that} \; \;  c^{-\la}_w c^\la_{\xi, v_\la} =
c^{-\la-\mu}_w c^{\la+\mu}_{\xi', v_{\la+\mu} }.
\end{equation}

For $y \in W_\g$ define the ideals
\begin{equation}
\label{id2}
Q(y)^\pm_w = \{c^{-\la}_w c^\la_{\xi, v_\la} \mid \xi \in V(\la)^*, \,
\xi \perp \UU_\pm T_y v_\la \}
\end{equation}
of $R^w_0$. We do not need to take span in the right hand side of 
\eqref{id2} because of \eqref{span}. The ideals \eqref{id2} are 
nontrivial if and only if $y \geq w$ in the plus case and 
$y \leq w$ in the minus case, see \cite{G}, and are completely
prime and $H$-invariant (with respect to the right action \eqref{Hact}).
Gorelik proved in \cite{G} that all $H$-invariant,  
prime ideals of $R^w_0$ are of the form $Q(y_-)^-_w + Q(y_+)^+_w$
for some $y_\pm \in W_\g$, $y_- \leq w \leq y_+$.

Recall that the quantum $R$-matrix associated to $w \in W_\g$ is 
given by
\begin{equation}
\RR^w = \prod_{j= k, \ldots, 1} \exp_{q_{i_j}}
\left( (1-q_{i_j}^{-2}) 
X^+_{\beta_j} \otimes X^-_{\beta_j} \right)
\label{Rw}
\end{equation}
where
\[
\exp_{q_i} (y)  = \sum_{n=0}^\infty q_i^{n(n+1)/2} 
\frac{y^n}{[n]_{q_i}!} \cdot
\]
and $q_i= q^{2/ \lcor \al_i, \al_i \rcor}$.
In \eqref{Rw} the terms are multiplied in the order $j = k, \ldots, 1$.
The $R$-matrix $\RR^w$ belongs to a certain completion \cite[\S 4.1.1]{L}
of $\UU^w_+ \otimes \UU^w_-$ and does not depend on the 
choice of a reduced decomposition of $w$. 

\bre{actions} The group $H$ acts on $\UU^w_-$ by $K_i . x = K_i x K_i^{-1}$, 
$x \in \UU^w_-$. The torus $(\KK^*)^{\times r}$ acts on $\UU^w_-$ 
by 
\[
(a_1, \ldots, a_r) \cdot x = 
(\prod_{i=1}^r a_i^{\lcor \la, \al_i\spcheck \rcor}) x, \quad
\forall x \in (\UU^w_-)_\la, \la \in Q_+,
\]
where $Q_+$ denotes the positive part 
$\Zset_{\geq 0} \al_1 + \ldots + \Zset_{\geq 0} \al_r$
of the root lattice of $\g$.
A subspace of $\UU^w_-$ is $H$-invariant if and only if 
it is $(\KK^*)^{\times r}$-invariant. In particular the 
set of $H$-primes and $(\KK^*)^{\times r}$-primes of $\UU^w_-$
are the same. Although $(\KK^*)^{\times r}$-invariance fits 
directly to the scheme in \cite{GL}, we will use the $H$-action 
since it is more intrinsic in term of the Hopf algebra
structure of $\UU_q(\g)$.
\ere

In \cite{Y} we proved:

\bth{hom1} \cite[Theorem 3.7]{Y} The map
\[
\phi_w \colon R_0^w \to \UU^w_-, \; \; 
\phi_w(c_w^{-\la} c^\la_{\xi, v_\la}) = 
(c^\la_{\xi, T_w v_\la} \otimes \id) (\RR^w), \; 
\la \in P_+, \xi \in V(\la)^* 
\]
is a (well defined) surjective algebra homomorhism.
It is $H$-equivariant with respect to the right action \eqref{Hact}
of $H$ on $R_0^w$. The kernel of $\phi_w$ is $Q(w)^+_w$.
\eth

This isomorphism is similar to one previously investigated 
by De Concini and Procesi in \cite{DP}.

Using Gorelik's description \cite{G} of the $H$-spectrum 
of $R^w_0$ leads to: 

\bth{Yold} \cite[Theorem 1.1]{Y} Fix $w \in W_\g$. For each $y \in W^{\leq w}$ define
\begin{multline}
\label{Iw}
I_w(y) = \phi_w (Q(y)_w^-+Q(w)_w^+) = \phi_w (Q(y)_w^-) \\
= \{ (c^{w,\la}_\eta \otimes \id)(\RR^w) \mid 
\la \in P_+, \eta \in (V_w(\la) \cap \UU_- T_y v_\la)^\perp \}.
\end{multline}
Then:

(a) $I_w(y)$ is an $H$-invariant completely prime ideal of $\UU^w_-$
and all $H$-invariant prime ideals of $\UU^w_-$ are
of this form.

(b) The correspondence $y \in W^{\leq w} \mapsto I_w(y)$ is an isomorphism
from the poset $W^{\leq w}$ to the poset of $H$-invariant 
prime ideals of $\UU^w_-$ ordered under inclusion; that is
$I_w(y) \subseteq I_w(y')$ for $y, y' \in W^{\leq w}$
if and only if $y \leq y'$.
\eth

Here and below for $w \in W_\g$ we denote by $W_\g^{\leq w}$ the set of all
$y \in W_\g$, $y \leq w$.
\sectionnew{Poisson structures on flag varieties}
\label{Poisson}

Denote by $\g_\Cset$ the complex simple Lie algebra 
corresponding to $\g$ and by $G_\Cset$
the connected, simply connected algebraic group with Lie algebra
$\g_\Cset$. Let $B^\pm_\Cset$ be a pair of opposite Borel 
subgroups. Set $T_\Cset = B^+_\Cset \cap B^-_\Cset$.
Denote by $\De_+$ the set of positive roots of $\g_\Cset$.
Fix two sets of root vectors $\{x^+_\al\}_{\al \in \De_+},$ 
$\{x^-_\al\}_{\al \in \De_+}$
with respect to $\Lie \, T_\Cset$ 
($x^\pm_\al \in \g_\Cset^{\pm \al}$) normalized by 
$\lcor x^+_\al, x^-_\al \rcor =1$, where $\lcor.,.\rcor$
denotes the 
nondegenerate invariant bilinear form on $\g_\Cset$
such that $\lcor \al, \al \rcor =2$ for a long root 
$\al$.

The standard Poisson structure \cite{GY} on the flag variety 
$G_\Cset/B^+_\Cset$ is defined
\begin{equation}
\label{pi}
\pi = \sum_{\al \in \De_+} \chi(x^+_\al) \wedge \chi(x^-_\al)
\end{equation}
where $\chi \colon \g_\Cset \to {\mathrm{Vect}}(G_\Cset/B_\Cset^+)$
denotes the infinitesimal action of $G_\Cset$ on $G_\Cset/B_\Cset^+$.
The action of the torus $T_\Cset$ on $(G_\Cset/B^+_\Cset, \pi)$ is 
Poisson.

The open Richardson varieties are the intersections of
opposite Schubert cells in the flag variety $G_\Cset/B_\Cset^+$
\begin{equation}
\label{Richardson}
R_{y_-, y_+} = B^-_\Cset y_- \cdot B^+_\Cset \cap 
B^+_\Cset y_+ \cdot B^+_\Cset \subset G_\Cset/B^+_\Cset,
\quad
y_\pm \in W.
\end{equation}
They are nonempty if and only if $y_- \leq y_+$ in the Bruhat order on $W_\g$.
In recent years Richardson varieties played an important role in
various algebro-geometric problems for flag varieties \cite{BL,KLS,BC} 
and the study of the totally nonnegative parts of flag varieties \cite{MR,RW}.

\bth{leaves}

(1) The $T_\Cset$-orbits of symplectic leaves of 
$(G_\Cset/B^+_\Cset, \pi)$ are precisely the open Richardson
varieties $R_{y_-, y_+}$, $y_\pm \in W, y_- \leq y_+$.
In particular, all open Richardson varieties are 
regular Poisson submanifolds of $(G_\Cset/B^+_\Cset, \pi)$.

(2) The codimension of a symplectic leaf 
in $R_{y_-, y_+}$ is
\[
\dim \ker (1 + y_+^{-1} y_-) =
\dim E_{-1}(y_+^{-1} y_-).
\]
\eth

\begin{proof} Part (1) is a special case of \cite[Theorem 0.4]{GY}, 
see \cite[Theorem 0.4]{BGY} in the case of 
Grassmannians. It also follows from \cite[Example 4.9]{EL}.

To deduce part (2), consider the double flag variety 
$G_\Cset/B_\Cset^+ \times G_\Cset/B_\Cset^-$ with the Poisson structure
\begin{multline*}
\pi^d = \sum_{\al \in \De_+} \big( 
\chi_1(x^+_\al) \wedge (\chi_1(x^-_\al) + \chi_2(x^-_\al) ) 
\\ - \chi_2(x^-_\al) \wedge (\chi_1(x^+_\al) + \chi_2(x^+_\al) )
\big) + \sum_j \chi_1(h_j) \wedge \chi_2(h_j). 
\end{multline*}
Here $\{ h_j\}$ denotes an orthonormal basis of 
$\Lie \, T_\Cset$ with respect to the restriction of $\lcor .,.\rcor$
and $\chi_{i} \colon \g_\Cset \to {\mathrm{Vect}}
(G_\Cset/B_\Cset^+ \times G_\Cset/B_\Cset^-)$,
$i=1,2$, denote the infinitesimal actions of $G_\Cset$ derived from the 
actions of $G_\Cset$ on the first and second factor of the Cartesian 
product.

It is easy to verify that the embedding 
of the single flag variety in the double flag variety 
$\eta \colon G_\Cset/B^+_\Cset \to G_\Cset/B_\Cset^+ \times G_\Cset/B_\Cset^-$
given by $\eta(g \cdot B^+_\Cset) = (g \cdot B^+_\Cset, B^-_\Cset)$ 
is Poisson with respect to $\pi$ and $\pi^d$. Obviously 
\begin{equation}
\eta( R_{y_-, y_+} ) = \left( (B^+_\Cset \times B^-_\Cset) \cdot 
(y_+ \cdot B^+_\Cset, B^-_\Cset) \right) 
\cap \left( \Delta(G_\Cset) \cdot
(y_- \cdot B^+_\Cset, B^-_\Cset) \right)
\label{emb}
\end{equation}
where $\Delta(G_\Cset)$ denotes the diagonal subgroup of 
$G_\Cset \times G_\Cset$.
Part (2) now follows from the fact \cite[Example 4.9]{EL} that 
the codimension of the symplectic leaves of the restriction 
of $\pi^d$ to the right hand side of \eqref{emb} 
is $\dim E_{-1}(y_+^{-1} y_-)$.
\end{proof}

The field $\Cset(R_{y_-, y_+})$ of rational functions 
on each open Richardson variety $R_{y_-, y_+}$ becomes
a Poisson field under the Poisson bracket induced 
from \eqref{pi}. Denote its center by $Z_\pi(\Cset(R_{y_-, y_+}))$.
The latter is called the field of rational 
Casimir functions on $(R_{y_-, y_+}, \pi)$. Each such function 
$f$ should be constant on the intersection of its domain 
with a generic symplectic leaf of $(R_{y_-, y_+}, \pi)$.
It is easy to see that the symplectic leaves of $\pi$ are 
smooth locally closed subvarieties of $G_\Cset/B^+_\Cset$,
e.g. by applying \cite[Theorem 1.9]{BGY}. Analogously to 
\cite[Proposition 4.2]{GeY} we obtain:

\ble{leq} For all $y_\pm \in W_\g$, $y_- \leq y_+$,
the transcendence degree of the field $Z_\pi(\Cset(R_{y_-, y_+}))$
is less than or equal to $\dim E_{-1}(y_+^{-1} y_-)$.
\ele

Another way to prove \leref{leq} is as follows.
Denote by $\pi^\sharp$ the bundle map 
$T^* R_{y_-, y_+} \to T R_{y_-, y_+}$ given by 
$\pi^\sharp_m (\eta_m) = (\eta_m \otimes \id) \pi_m$ 
for $m \in R_{y_-, y_+}$, $\eta_m \in T^*_m R_{y_-, y_+}$.
Then 
\[
Z_\pi(\Cset(R_{y_-, y_+})) = \{ 
f \in \Cset(R_{y_-, y_+}) \mid \\
\pi^\sharp(\eta) (f) = 0, \; \forall 
\eta \in \Ga( R_{y_-, y_+},  T^*_m R_{y_-, y_+}) \}.
\]
Since the codimension of all symplectic leaves of 
$(R_{y_-, y_+}, \pi)$ is $\dim E_{-1}(y_+^{-1} y_-)$
there exist $k= \dim R_{y_-, y_+} - \dim E_{-1}(y_+^{-1} y_-)$
generically linearly independent regular vector fields 
$X_1, \ldots, X_k$ on $R_{y_-, y_+}$ such that 
\[
X_j f = 0, \quad \forall f \in Z_\pi(\Cset(R_{y_-, y_+})), 
j=1, \ldots, k.
\]
This implies that the transcendence degree of 
$Z_\pi(\Cset(R_{y_-, y_+}))$ is less than or equal to 
$\dim R_{y_-, y_+} -k = \dim E_{-1}(y_+^{-1} y_-)$.

In \thref{main} we will prove an inverse inequality which 
will show that the transcendence degree of $Z_\pi(\Cset(R_{y_-, y_+}))$
is equal to $\dim E_{-1}(y_+^{-1} y_-)$.
This can be also obtained directly by constructing 
enough functions in the Poisson center without going to 
quantized algebras of functions.
\sectionnew{Dimensions of the Goodearl--Letzter strata}
\label{dimensions}
Fix $w \in W_\g$ and $y \in W_\g^{\leq w}$. By Theorems \ref{thom1}
and \ref{tYold}
\begin{equation}
\UU^w_-/I_w(y) \cong R_0^w/(Q(y)_w^- + Q(w)_w^+).
\label{isom}
\end{equation}
Since $Q(y)_w^- + Q(w)_w^+$ is a completely prime ideal, the 
quotient rings in \eqref{isom} are domains.
Denote by $L_{y,w}$ the localization of 
$R_0^w/(Q(y)_w^- + Q(w)_w^+)$ by all nonzero 
homogeneous elements with respect to the right 
$H$-action \eqref{Hact}.
The Goodearl--Letzter results \cite[Theorem 6.6]{GL}
and the Brown--Goodearl result on strong rationality 
\cite[Theorem II.6.4]{BG} imply that the 
center $Z(L_{y,w})$ of $L_{y,w}$
is a Laurent polynomial ring over $\KK$.
Denote by $n_{y,w}$ 
the number of independent variables in this ring. 
Then \cite[Theorem 6.6]{GL}
implies that $\Spec_{I_w(y)} \UU^w_-$ is homeomorphic to 
the spectrum of a Laurent polynomial ring in 
$n_{y,w}$ variables. In this section 
we prove that $n_{y,w} = \dim E_{-1}(w^{-1} y)$.
  
For $\la \in P_+$ set
\begin{equation}
a_\la = c_w^{-\la} c_y^\la,
\label{a-la}
\end{equation}
recall \eqref{cwla}.
Applying the standard $R$-matrix commutation relations 
\cite[Theorem I.8.16]{J} in $R_q[G]$ gives
\begin{multline}
a_\la \big( c_w^{-\mu} c^{\mu}_{\xi, v_\mu} \big) 
- q^{- \lcor(y+w)\la, \nu + w\mu \rcor} \big( c_w^{-\mu} c^{\mu}_{\xi, v_\mu} \big)
a_\la \in Q(y)_w^- + Q(w)_w^+, \\
\forall \xi \in V(\mu)^*_\nu, \mu \in P_+, \nu \in P. 
\label{comm}
\end{multline}
Denote the image of $a_\la$ in 
$R_0^w/(Q(y)_w^- + Q(w)_w^+)$ by $\ol{a}_\la$.
By \eqref{comm} all $\ol{a}_\la$
are normal elements. They are all nonzero
since for $w_1, w_2 \in W_\g$, $\la \in P_+$, one has
$T_{w_1} v_\la \in \UU_+ T_{w_2} v_\la$ 
if and only if $w_1 \geq w_2$, see \cite[Proposition 4.4.5]{J}.
 
Recall that $L_{y,w}$ denotes the localization of 
$R_0^w/(Q(y)_w^- + Q(w)_w^+)$ by the set of 
all nonzero homogeneous elements.
Represent $\la \in P$ as $\la = \la_+ - \la_-$, 
for $\la_{\pm} \in P_+$ with non-intersecting support and
set
\[
\ol{a}_\la = (\ol{a}_{\la_-})^{-1} \ol{a}_{\la_+} \in L_{y,w}.
\]

Then $\ol{a}_\la$ are normal elements of $L_{y,w}$ for all 
$\la \in P$ and 
\begin{multline}
\ol{a}_\la \big( c_w^{-\mu} c^{\mu}_{\xi, v_\mu} + Q(y)_w^- + Q(w)_w^+\big) 
= \\
q^{- \lcor(y+w)\la, \nu + w\mu \rcor} 
\big( c_w^{-\mu} c^{\mu}_{\xi, v_\mu} + Q(y)_w^- + Q(w)_w^+ \big)
\ol{a}_\la 
\label{comm2}
\end{multline}
for all
$\xi \in V(\mu)^*_\nu, \mu \in P_+, \nu \in P$. 

Denote 
\[
P_{y,w} = \{ \la \in P \mid (y+w) \la = 0 \}.
\]
Obviously $P_{y,w}$ is a lattice of rank 
$\dim E_{-1}(w^{-1}y) = \dim \ker(y+w)$. 
Denote by $\la_1, \ldots, \la_{ \dim E_{-1}(w^{-1} y)  }$ 
a basis of $P_{y,w}$. Set
\[
Z_{y,w} = \Span \{ \ol{a}_\la \mid \la \in P_{y,w} \}.
\]
\ble{1} For all $y \leq w$, $Z_{y,w}$ is a subring 
of $Z(L_{y,w})$ which is isomorphic to a Laurent 
polynomial ring in $\dim E_{-1} (w^{-1} y)$ variables
over $\KK$. In particular, $n_{y,w} \geq \dim E_{-1} (w^{-1} y)$.
\ele
\begin{proof}
The fact that $Z_{y,w} \subseteq Z(L_{y,w})$ follows 
from \eqref{comm}. Applying \cite[eq. (3.14)]{Y}, one gets that for all 
$\la, \mu \in P$, $\ol{a}_\la \ol{a}_\mu$ is a nonzero 
scalar multiple of $\ol{a}_{\la+\mu}$. In particular 
$Z_{y,w}$ is a ring generated by 
$\left( \ol{a}_{\la_i} \right)^{\pm 1}$, 
$i =1, \ldots, \dim E_{-1} (w^{-1} y)$. It is isomorphic to 
a Laurent polynomial ring in $\dim E_{-1} (w^{-1} y)$ 
variables since the elements 
$\ol{a}_\la$ are linearly independent for different 
$\lambda$'s. This is so because the elements $\ol{a}_\la$ 
belong to different homogeneous components of $L_{y,w}$
with respect to right action \eqref{Hact} of $H$:
\[
\ol{a}_\la \lha K_i = 
q^{\lcor (y-w)\la, \al_i \rcor } \ol{a}_\la = 
q^{2 \lcor y \la, \al_i \rcor } \ol{a}_\la 
\]
for all $\la \in P_{y,w}$.
\end{proof}

The main result of the paper is:

\bth{main} Assume that the base field $\KK$ has
characteristic $0$ and that the deformation parameter
$q$ is transcendental over $\Qset$.

For all $w \in W_\g$, $y \in W^{\leq w}$ the
Goodearl--Letzter stratum $\Spec_{I_w(y)} \UU^w_-$ 
is homeomorphic to the spectrum of a Laurent 
polynomial ring over $\KK$ in $\dim E_{-1}(w^{-1} y)$ 
variables and the transcendence degree of 
$Z_\pi(\Cset(R_{y, w}))$ is equal 
to $\dim E_{-1}(w^{-1} y)$.
\eth

As we noted earlier, \cite[Theorem II.6.4]{BL} 
implies that the stratum $\Spec_{I_w(y)} \UU^w_-$
is homeomorphic to the spectrum of a Laurent
polynomial ring over $\KK$.
It is sufficient to prove the first part of 
\thref{main} for one choice of the base field $\KK$ 
since all algebras $R_q[G], R_0^w, Q(y)_w^- + Q(w)_w^+, 
L_{y,w}, Z(L_{y,w})$ behave appropriately under 
extensions of the base field from $\Qset(q)$ to $\KK$.
We will prove \thref{main} for $\KK = \Cset(q)$
and from now on we will assume that $\KK= \Cset(q)$.

Before we proceed with the proof, we will recall some 
facts about integral forms. Set $\AA = \Cset[q, q^{-1}]$.
Denote by $\UU^{\res}_\AA$ the restricted integral form 
of $\UU_q(\g)$ over $\AA$, see \cite[Chapter 9]{CP}
for detail. (Usually integral forms of quantized
universal enveloping algebras are defined over 
$\Zset[q,q^{-1}]$, but this will not be needed for 
our purposes.) Denote 
$V(\la)^{\res}_\AA = \UU^{\res}_\AA v_\la \subset V(\la)$,
cf. \cite[\S 10.1]{CP}. Note that 
$\UU^{\res}_\AA$ and $V(\la)^{\res}_\AA$ are stable under 
the action of the braid group $\BB_\g$ for 
all $\la \in P_+$. In particular, $V(\la)^{\res}_\AA$
and $(V(\la)^*)^{\res}_\AA$
are also generated by lowest weight vectors: 
$V(\la)^{\res}_\AA= \UU^{\res}_\AA T_{w_\ci} v_\la$, 
$(V(\la)^*)^{\res}_\AA = \UU^{\res}_\AA \xi_\la$. 
Denote by $R^+_\AA$ the $\AA$-subalegbra of $R^+$ 
consisting of all sums of 
elements of the form $c^\la_{\xi, v_\la}$, 
for $\la \in P_+$, 
$\xi \in (V(\la)^*)^{\res}_\AA = \UU^{\res}_\AA \xi_\la$.
Note that $c_w \subset R^w_\AA$ and
denote $R^w_\AA = R^+_\AA[c_w^{-1}]$. The group 
$H$ acts on $R^w_\AA$ on the left and right 
by \eqref{Hact}. The invariant subalgebra with 
respect to the left action will be denoted by 
$(R^w_0)_\AA$. Clearly,
\begin{equation}
\label{Rwaa}
(R^w_0)_\AA = \{ c_w^{-\la} c^\la_{\xi, v_\la} \mid 
\la \in P_+, \xi \in (V(\la)^*)^{\res}_\AA\}.
\end{equation}
Analogously to \eqref{Rw0} one does not need to take
a sum in the right hand side of \eqref{Rwaa}.
Similarly to \eqref{id2}, for $u \in W_\g$ define the ideals
\[
(Q(u)^\pm_w)_\AA = \{c^{-\la}_w c^\la_{\xi, v_\la} \mid 
\xi \in (V(\la)^*)^{\res}_\AA, \, \xi \perp \UU_\pm T_y v_\la \}.
\]
It is clear that the natural embeddings
\[
R^+_\AA \otimes_\AA \Cset(q) \to R^+, \; \; 
(R^w_0)_\AA \otimes_\AA \Qset(q) \to R^w_0, \; \; 
(Q(u)^\pm_w)_\AA \otimes_\AA \Cset(q) \to Q(u)^\pm_w, 
\]
and 
\[
[(R^w_0)_\AA/((Q(y)^-_w)_\AA + (Q(w)^+_w)_\AA) ] \otimes_\AA 
\Cset(q) \to R_0^w/(Q(y)_w^- + Q(w)_w^+)
\]
are isomorphisms. Denote the localization $\AA_1 = (\AA)_{(q-1)}$
of $\AA$ at $(q-1)$. Let $(L_{y,w})_{\AA_1}$ denote 
the localization of $(R^w_0)_\AA/((Q(y)^-_w)_\AA + (Q(w)^+_w)_\AA) )$
by all homogeneous elements which do not belong to
\[ 
(q-1) (R^w_0)_\AA/((Q(y)^-_w)_\AA + (Q(w)^+_w)_\AA ).
\]
Then $(L_{y,w})_{\AA_1}$ has a canonical structure of 
an $\AA_1$ algebra and 
\begin{equation}
\label{inc}
(L_{y,w})_{\AA_1} \otimes_{\AA_1} \Cset(q) \cong L_{y,z} 
\; \; \mbox{and} \; \; 
Z((L_{y,w})_{\AA_1}) \otimes_{\AA_1} \Cset(q) \cong Z(L_{y,z}).
\end{equation}

Next, we recall some general facts about specializations. 
Assume that $C$ is an $\AA$ algebra. One calls the $\Cset$-algebra
\begin{equation}
\ol{C}=C/(q-1)C \cong C \otimes_\AA \Cset
\end{equation}
the specialization 
of $C$ at $1$. The tensor product is defined via the homomorphism
$\kappa \colon \AA \to \Cset$, $\kappa(q)=1$. If $\ol{C}$ is a commutative 
algebra then it has a canonical Poisson algebra structure defined as follows.
Denote the quotient map $\eta \colon C \to C/(q-1)C = \ol{C}$
For $a, b \in \ol{C}$ let $a', b'$ be two preimages of $a,b$ 
under $\eta$. Define  
\begin{equation}
\label{Po}
\{a,b\} = \eta( (a' b' - b' a')/(q-1) ). 
\end{equation}
It is well known that the above bracket does not depend
on the choice of $a', b'$ and that it turns $\ol{C}$ into a Poisson 
algebra. Obviously
\begin{equation}
\label{eta}
\eta(Z(C)) \subseteq Z_{\{.,.\}}(\ol{C}) 
\end{equation}
where the left hand side denotes the center of the 
Poisson algebra $(C, \{.,.\})$. 

The same construction defines the 
specialization $\ol{C}$ of an $\AA_1$ algebra $C$ at $1$. 
Provided that $C/(q-1)C$ is commutative, one defines 
a Poisson algebra structure on $\ol{C}$ by \eqref{Po}.

It is well known that 
\begin{equation}
\label{indR+}
\ol{R^+} \cong \Cset[G_\Cset/U^+_\Cset]
\end{equation}
where $U^+_\Cset$ denotes the unipotent radical of $B^+_\Cset$.
This is implicit in \cite[\S 9.1.6]{J0} and appears in the 
proof in \cite[Th\'eor\`eme 3]{J2}.
Eq. (4.9) and Theorem 4.6 (both) from \cite{Y} imply at once that 
\begin{equation}
\label{Psp}
\ol{R_0^w} \cong \Cset[w B^-_\Cset \cdot B^+_\Cset]
\end{equation}
(where $w B^-_\Cset \cdot B^+_\Cset \subset G_\Cset/B^+_\Cset$ 
is the translated one Schubert cell of $G_\Cset/B^+_\Cset$, 
see \cite[\S 4.2]{Y}) and
\begin{equation}
\label{1isom}
\ol{(R_0^w/(Q(y)_w^- + Q(w)_w^+))_\AA} \cong \Cset[R_{y,w}].
\end{equation}
Note that in the setting of \cite[Theorem 4.6]{Y}, 
the variety $S_w(y)$ is isomorphic to $R_{y,w}$ via 
the map (4.10). We claim that the induced Poisson structure 
on $\Cset[R_{y,w}]$ is exactly the one coming from $- \pi$,
recall \eqref{pi}.
The above argument proves that it is sufficient to show 
that the induced from \eqref{Psp} Poisson structure on 
$\Cset[w B^-_\Cset \cdot B^+_\Cset]$ is equal to 
the one coming from $- \pi$ which follows from 
the following lemma.
\ble{ind} The induced from \eqref{indR+} Poisson structure on
$\Cset[G/ U^+_\Cset]$ is equal to the one coming from
following bivector field on $G/ U^+_\Cset$
\[
- \sum_{\al \in \De_+} \chi(x_\al^+) \wedge \chi(x_\al^-),
\]
where $\chi \colon \g_\Cset \to G_\Cset /U^+_\Cset$ is the 
infinitesimal action of $G_\Cset$ on $G_\Cset/U^+_\Cset$. 
\ele
\begin{proof} Denote by $w_0$ the longest element of $W_\g$. 
The algebras $\UU_\pm$ are $Q$ graded by 
$\deg X_i^\pm = \pm \al_i$. Denote by $(\UU_\pm)_{\pm \ga}$, 
$\ga \in Q_+$ their graded components and by 
$m(\ga) = \dim (\UU_+)_\ga= \dim (\UU_-)_{-\ga}$ 
their dimensions. For each $\ga \in Q_+$ 
fix a pair of dual bases
$\{u_{\ga, i} \}_{i=1}^{m(\ga)}$ and
$\{u_{-\ga, i} \}_{i=1}^{m(\ga)}$
of $(\UU_+)_\ga$ and $(\UU_-)_{-\ga}$
with respect to the Rosso--Tanisaki form.
Then
\begin{equation}
\label{Rm}
\RR^{w_0} = 1+ \sum_{\ga \in Q_+, \ga \neq 0}
\sum_{i=1}^{m(\ga)} u_{\ga, i} \otimes u_{-\ga, i}
\end{equation}
and we have the $R$-matrix commutation relation in 
$R^+$:
\begin{multline}
\label{commr}
c_{\xi_1, v_{\la_1}}^{\la_1} c_{\xi_2, v_{\la_2}}^{\la_2} =
q^{ \lcor \nu_1, \nu_2 \rcor - \lcor \la_1, \la_2 \rcor}
\Big(
c_{\xi_2, v_{\la_2}}^{\la_2} c_{\xi_1, v_{\la_1}}^{\la_1}
\\
+ \sum_{\ga \in Q_+, \ga \neq 0}
\sum_{i=1}^{m(\ga)}
c_{S^{-1}(u_{\ga, i})\xi_2, v_{\la_2}}^{\la_2}
c_{S^{-1}(u_{-\ga, i}) \xi_1, v_{\la_1}}^{\la_1}
\Big), 
\end{multline}
for $\la_i \in P_+$, $\nu_i \in P$,
$\xi_i \in V(\la_i)^*_{\nu_i}$, $i=1,2$.
Eq. \eqref{Rw} implies:
\[
\ol{(R^{w_0} -1)/(q-1)} = 2 \sum_{\al \in \De_+} 
x_\al^+ \otimes x_\al^-.  
\] 
Therefore in terms of the notation from \eqref{commr} 
\begin{multline*}
\Big\{ {\ol{c_{\xi_1, v_{\la_1}}^{\la_1}}}, 
{\ol{c_{\xi_2, v_{\la_2}}^{\la_2}}} \Big \}
= (\lcor \nu_1, \nu_2 \rcor - \lcor \la_1, \la_2 \rcor ) 
{\ol{c_{\xi_1, v_{\la_1}}^{\la_1}}} \;
{\ol{c_{\xi_2, v_{\la_2}}^{\la_2}}} 
\\ + 2 \sum_{\al \in \De_+} 
\Big( \chi(x_\al^-) {\ol{c_{\xi_1, v_{\la_1}}^{\la_1}}} \Big) 
\Big( \chi(x_\al^+) {\ol{c_{\xi_2, v_{\la_2}}^{\la_2}}} \Big)
\end{multline*}
and
\begin{multline*}
\Big\{ {\ol{c_{\xi_1, v_{\la_1}}^{\la_1}}},
{\ol{c_{\xi_2, v_{\la_2}}^{\la_2}}} \Big \}
= \frac{1}{2} \Big\{ {\ol{c_{\xi_1, v_{\la_1}}^{\la_1}}},
{\ol{c_{\xi_2, v_{\la_2}}^{\la_2}}} \Big \} - 
\frac{1}{2} 
\Big\{ {\ol{c_{\xi_2, v_{\la_2}}^{\la_2}}},
{\ol{c_{\xi_1, v_{\la_1}}^{\la_1}}} \Big \}
\\ = 
- \sum_{\al \in \De_+}
\Big( \chi(x_\al^+) {\ol{c_{\xi_1, v_{\la_1}}^{\la_1}}} \Big)
\Big( \chi(x_\al^-) {\ol{c_{\xi_2, v_{\la_2}}^{\la_2}}} \Big)
+
\sum_{\al \in \De_+}
\Big( \chi(x_\al^-) {\ol{c_{\xi_1, v_{\la_1}}^{\la_1}}} \Big)
\Big( \chi(x_\al^+) {\ol{c_{\xi_2, v_{\la_2}}^{\la_2}}} \Big)
\end{multline*}

\end{proof}
The above proves that $\ol{(L_{y,w})_{\AA_1}}$ is isomorphic to a 
subring of $\Cset(R_{y,w})$ which contains 
$\Cset[R_{y,w}]$ and the induced Poisson bracket 
on $\ol{(L_{y,w})_{\AA_1}}$ coincides with the 
one coming from $- \pi$.
\\ \hfill \\
\noindent
{\em{Proof of \thref{main}}.} Assume that the first statement is not 
true. \leref{1} implies that $n_{y,w} \geq E_{-1} (w^{-1} y) +1$.
From \eqref{inc} one obtains that the Krull dimension of 
$Z((L_{y,w})_{\AA_1})$ is at least $\dim E_{-1} (w^{-1} y) +2$. 
On the other hand \eqref{eta} and the discussion 
before the proof of the theorem imply that 
\[
Z((L_{y,w})_{\AA_1})/ [(q-1) Z((L_{y,w})_{\AA_1})]
\] 
is isomorphic to a subring of $Z_\pi(\Cset(R_{y_-, y_+}))$. 
It follows from \leref{leq} that the Krull dimension of 
$Z((L_{y,w})_{\AA_1})$ is at most $\dim E_{-1} (w^{-1} y) +1$
which is a contradiction. This proves the first part of the theorem.

The second part of the theorem follows from the first. From 
\leref{leq} we know that the transcendence degree of 
$Z_\pi(\Cset(R_{y, w}))$
is less than or equal to $\dim E_{-1}(w^{-1} y)$.
If the inequality is strict, then working the above argument 
backwards leads to $n_{y,w} <  \dim E_{-1}(w^{-1} y)$ which is a 
contradiction. This completes the proof of the theorem.
\qed 

\bre{com} Due to Joseph \cite{J0} and Hodges--Levasseur--Toro \cite{HLT}
the $H$-prime ideals of the quantized coordinate 
ring $\OO_q(G)$ are parametrized by pairs $(y,w) \in W_\g \times W_\g$ 
and the corresponding stratum of $\Spec \OO_q(G)$ has dimension
equal to 
$\dim E_1(w^{-1} y)$. On the Poisson side the 
$T_\Cset$-orbits of symplectic 
leaves of the Poisson--Lie group $G_\Cset$ equipped with the  
standard Poisson structure are the double Bruhat cells 
$G^{y,w}_\Cset= B^-_\Cset y B^-_\Cset 
\cap B^+_\Cset w B^+_\Cset$ and the codimension 
of a symplectic leaf of $G^{y,w}_\Cset$ in $G^{y,w}_\Cset$ is equal
\cite[Theorem A.2.1 and Proposition A.2.2]{HL} 
to $\dim E_{1}(w^{-1} y)$. It is not hard to show 
that the tanscendence degree of the center of the Poisson 
field $\Cset(G^{y,w}_\Cset)$ is equal to $\dim E_{1}(w^{-1} y)$.

It would be very interesting to understand the relation between 
these facts and the dimension formulas in \thref{main}.
(Note the difference of $\pm 1$ eigenspaces.)
We believe 
that there exists a ring theoretic counterpart of the construction 
of weak splitting of surjective Poisson submersions in 
\cite[Sect. 3]{GY} which relates the $H$-stratifications 
of $\Spec \OO_q(G)$ and $\Spec \UU^w_-$, and eventually 
the dimension formulas.  
\ere

\end{document}